\documentclass{article}
\usepackage[utf8]{inputenc}
\usepackage[english]{babel}
\usepackage{graphicx}
\usepackage{amssymb}
\usepackage{amsmath}
\usepackage{commath}
\usepackage[hidelinks]{hyperref}
\usepackage{color}
\usepackage{wrapfig}
\usepackage{booktabs}
\usepackage{graphicx,wrapfig,lipsum}
\graphicspath{ {./images/} }
\usepackage[a4paper, total={6in, 8in}]{geometry}

\title{A new coupled model for laser-induced thermal ablation accounting for tissue water vaporization}

\author{Federico Herrero-Hervás$^a$ \and \hspace{1 cm} Angiolo Farina$^b$}
\date{ }

\begin{document}

\maketitle

\maketitle
$^a$ Department of Mathematical Analysis and Applied Mathematics, Instituto de Matemática 

\hspace{0.2 cm} Interdisciplinar, Universidad Complutense de Madrid.

$^b$ Dipartimento di Matematica e Informatica `Ulisse Dini', Università degli Studi di Firenze.

\begin{center}
\large \textbf{Abstract}     
\end{center}
A new model is considered to describe the evolution of tissue temperature and water concentration during a thermal ablation process. 
The aim of this study is to expand the results obtained in \cite{LITT}, where the tissue water concentration is considered as a known function, experimentally obtained.  
\\\\
The model consists of two coupled equations: a parabolic partial differential equation for the temperature $\theta$ and an ordinary differential equation for the water content $w$.
\\\\
The equation for $\theta$ accounts for temperature diffusion, blood perfusion, the heat source provided by the laser applicator and the vaporization effects which act when the temperature approaches 100ºC. The equation for $w$ accounts for the loss of tissue water that occurs during vaporization.
\\\\
The model equations are solved in a cylindrical domain around the laser applicator and the results are compared with the existing experimental data for porcine liver tissue from \cite{4}.

\section{Introduction}
Laser-induced thermal ablation is a method used for tumor destruction based on heat ablation of the tumorous tissue. The main feature of this therapy relies on the placement of an applicator equipped with a laser radiating tip into the tumorous tissue. The energy source produced by the irradiation of the laser results in a temperature increase, with the aim of destructing the tumorous tissue. 

The foremost underlying processes are the laser irradiation emitted from the applicator, the thermal diffusion along the tissue, blood perfusion when the treatment is applied to an \textit{in-vivo} organ and the vaporization of tissue water. It is precisely tissue water vaporization the part that is most difficult to represent. Initial studies such as \cite{Fasano, 3} did not include it for their models, however for temperatures reaching 100 ºC, the effect becomes more relevant. In more recent papers like \cite{LITT}, it is addressed in different ways, mainly taking tissue water concentration as a known function, experimentally fitted.

Instead of considering the fitted function representing tissue water concentration, we propose a model which treats the water vaporization effects as a coupled term modulating the temperature evolution. The model hence consists of two variables, the tissue temperature and water concentration, that are obtained as the solutions to a coupled system of differential equations. 

We model the process as a parabolic-ordinary system of differential equations, estimating its parameters with the literature data. The system is then numerically solved and compared with experimental data for \textit{ex-vivo} porcine livers.

The paper is structured as follows: the main mathematical model for the ablation treatment is described on Section \ref{s1}. We start from Penne's bioheat equation and add the relevant terms, especially the modulation effect caused by tissue water vaporization. In Section \ref{param} the specifics of the source and vaporization terms are detailed, as well as the main parameters for porcine liver tissue. The results of the numerical simulation of the process are gathered on Section \ref{s3}, with a comparison to the already existing experimental data. In Section \ref{s4}, some final remarks are discussed.

\section{Modeling the thermal ablation procedure} \label{s1}
The main purpose of our model is to obtain two coupled equations that describe the evolution of the tissue temperature and water concentration in an organ subject to laser irradiation over time. Let $\Omega \subset \mathbb{R}^3$ be the tissue domain. We denote by $\theta(\Vec{x},t)$ and $w(\Vec{x},t)$ the temperature and tissue water content (i.e. the water fraction per unit mass of tissue) on $\Vec{x} \in \Omega$, $t>0$, respectively.

The equation for $\theta$ includes thermal diffusion, while also accounts for the source term due to the laser irradiation, blood perfusion and the effects of tissue water vaporization. These elements are incorporated to Penne's bioheat equation (see for instance \cite{Hemomath}). As for $w$, an evolution equation is required to model the loss of water as a result of the increasing tissue temperature. 

\subsection{Bioheat equation}

As a first approximation to the model, we start by neglecting the vaporization effects, that are only meaningful at large temperatures approaching 100ºC. For mild temperatures, the resulting equations are of the form
\begin{equation} \label{1}
\begin{cases}
    \rho _{ts}c_{ts}\dfrac{\partial \theta }{\partial t}=k_{ts}\Delta \theta
-\rho _{bl}c_{bl}\omega \left( \theta -\theta _{bl}\right) + S(\Vec{x},t), \\
\\
\dfrac{\partial w}{\partial t}=-\Gamma \left( \theta \right) w,
\end{cases} \quad \text{for } \Vec{x} \in \Omega, ~t>0,
\end{equation}
where $\rho _{\alpha }$, $\alpha =ts,~bl$, is the mass density of the tissue, blood and water, $c_{\beta }$, $\beta =ts,~bl$, is the specific heat of the tissue and blood, $\omega $ is the blood perfusion rate, $\theta _{bl} $ is the blood temperature, $k_{ts}$ is the tissue heat diffusivity coefficient, $S$ is the source function modeling laser irradiation and $\Gamma $ is the vaporization rate depending on the local temperature.

\begin{figure}[htp]
    \centering
    \includegraphics[scale=0.7]{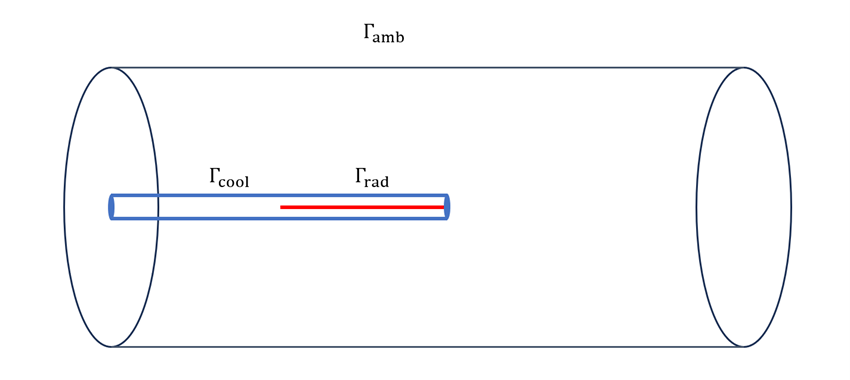}
    \caption{Geometry of the domain $\Omega$ and its boundary.}
    \label{fig:f1}
\end{figure}

Equations \eqref{1} require boundary conditions for $\theta$ and the initial tissue temperature and water concentration distributions. To set these conditions, we impose a geometrical restriction on $\Omega$, considering a large enough cylindrical domain, with the laser applicator on its axis, as in Figure \ref{fig:f1}.
\\The applicator is surrounded by a cooling flow of water to prevent burns in the tissue. The irradiation is emitted from the last end of the applicator, where a radiating laser fiber is placed. These are schematically represented in blue and red in Figure \ref{fig:f1}. Hence, the boundary of $\Omega$ is made up of the outer face of the cylinder, denoted by $\Gamma_{\text{amb}}$, where a heat exchange with the ambient temperature has to be considered; and the boundary corresponding to the applicator. This is subsequently split into two parts, one denoted $\Gamma_{\text{cool}}$, where only the coolant is present, and the second one, $\Gamma_{\text{rad}}$ in the part around the radiating tip of the applicator.

\begin{figure}[htp]
    \centering
    \includegraphics[scale=0.4]{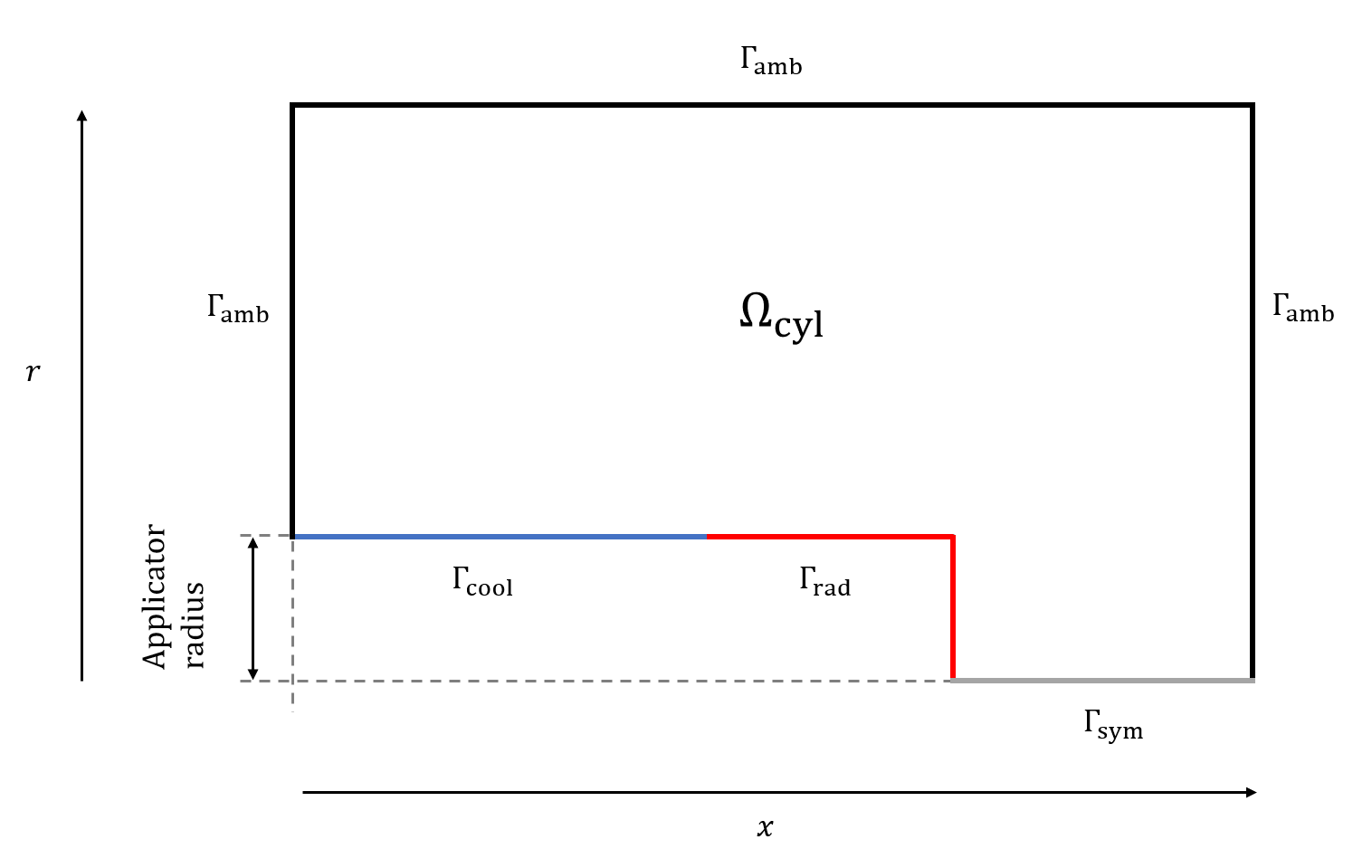}
    \caption{Two dimensional domain after the considered cylindrical symmetry.}
    \label{fig:f2}
\end{figure}

We further assume axial symmetry of with respect to the axis of the cylinder and thus only look for solutions depending on $r$ (the radial distance), $x$ (the longitudinal coordinate), and time $t$. This way, the two-dimensional domain considered, denoted $\Omega_{\text{cyl}}$ is as depicted in Figure \ref{fig:f2}.

The boundary of $\Omega_{\text{cyl}}$, apart from the previous $\Gamma_{\text{amb}}$, $\Gamma_{\text{cool}}$ and $\Gamma_{\text{rad}}$, also consists of the continuation of the cylinder's axis through the tissue after the applicator ends, denoted by $\Gamma_{\text{sym}}$ (see again Figure \ref{fig:f2}).

The boundary conditions on $\Gamma_{\text{amb}}$, $\Gamma_{\text{cool}}$ and $\Gamma_{\text{rad}}$ must reflect the heat exchange between the outer temperature (ambient temperature for $\Gamma_{\text{amb}}$ and coolant temperature for $\Gamma_{\text{cool}}$ and $\Gamma_{\text{rad}}$) and $\theta$, whereas on $\Gamma_{\text{sym}}$ a symmetry condition has to be set for the axial symmetry to hold. 

We thus have
\begin{equation} \label{2}
    \begin{cases}
        k_{ts} \dfrac{\partial \theta}{\partial \Vec{n}} = \alpha_{\text{amb}}(\theta_{\text{amb}} - \theta) \quad & \text{on } \Gamma_{\text{amb}}, \\\\
        k_{ts} \dfrac{\partial \theta}{\partial \Vec{n}} = \alpha_{\text{cool}}(\theta_{\text{cool}} - \theta) \quad & \text{on } \Gamma_{\text{rad}} \cup \Gamma_{\text{cool}}, \\\\
        \dfrac{\partial \theta}{\partial \Vec{n}} = 0 \quad & \text{on } \Gamma_{\text{sym}},
    \end{cases}
\end{equation}
where $\theta_{\text{amb}}$ and $\theta_{\text{cool}}$ are the known ambient and cooling temperatures, while $\alpha_{\text{amb}}$ and $\alpha_{\text{cool}}$ are the heat transfer coefficients for the tissue and the cooled part of the applicator, respectively, and $\Vec{n}$ is the normal outward vector.

The equation for $w$ does not require boundary conditions, as no spatial derivatives are involved. Both however need initial values, $\theta(r,x,0)$ and $w(r,x,0)$. To match the \textit{ex-vivo} experiments illustrated in \cite{4}, the initial temperature distribution $\theta(r,x,0)$ is set to $\theta_{\text{amb}}$, while the initial tissue water concentration is discussed in more detail on Section \ref{s1.3}.

\subsection{Model for the vaporization effects}

After considering the initial model of equations \eqref{1}, the effects of vaporization have to be taken into account. The results obtained in \cite{4} show the need to include some modulation effects in the equation for $\theta$, as the temperature cannot constantly increase while the tissue water is vaporizing. The model proposed in \cite{LITT} incorporated vaporization through two different approaches: an effective specific heat model, where tissue water density was taken as a known function of $\theta$; and an enthalpy model, with a second equation for the enthalpy.

Our approach however aims to modulate the temperature growth through the tissue water concentration $w$, which, in turn, evolves in time according to equation \eqref{1}$_2$. This way, we propose a modification of system \eqref{1} as follows
\begin{equation} \label{3}
\begin{cases}
    \rho _{ts}c_{ts}\dfrac{\partial \theta }{\partial t}=k_{ts}\Delta \theta
-\rho _{bl}c_{bl}\omega \left( \theta -\theta _{bl}\right) + S(\Vec{x},t) \cdot f(\theta,w), \\
\\
\dfrac{\partial w}{\partial t}=-\Gamma \left( \theta \right) w,
\end{cases} \quad \text{for } \Vec{x} \in \Omega, ~t>0,
\end{equation}
where we have introduced a modulating function $f$ depending both on $\theta$ and $w$ that dampens the temperature growth caused by the source term $S$. The expression for $f$ must model the flattening of the temperature growth as it approaches 100ºC, where the thermal energy goes into water vaporization instead of causing temperature increase. Once the tissue water has vaporized, the energy is again used for the temperature increase. Thus, our proposed function is of the form
\begin{equation} \label{4}
    f(\theta,w) = \begin{cases}
       \left(\dfrac{\theta_{\text{vap}_2} - \theta}{\theta_{\text{vap}_2}- \theta_{\text{vap}_1}} \right)^\gamma \quad & \text{if} ~\theta_{\text{vap}_1} \leq \theta \leq \theta_{\text{vap}_2} ~\text{and } w > w_\text{vap}, \\\\
       1 \quad & \text{otherwise},
    \end{cases}
\end{equation}
where $\theta_{\text{vap}_1} < \theta_{\text{vap}_2}$ are the two threshold values in between which vaporization takes place, the exponent $\gamma>0$ is a parameter and $0 < w_\text{vap} \ll 1$ is another threshold value for the tissue water concentration from which we can consider the remaining amount of water still to be vaporized as negligible.

\begin{figure}[htp]
    \centering
    \includegraphics[scale=0.45]{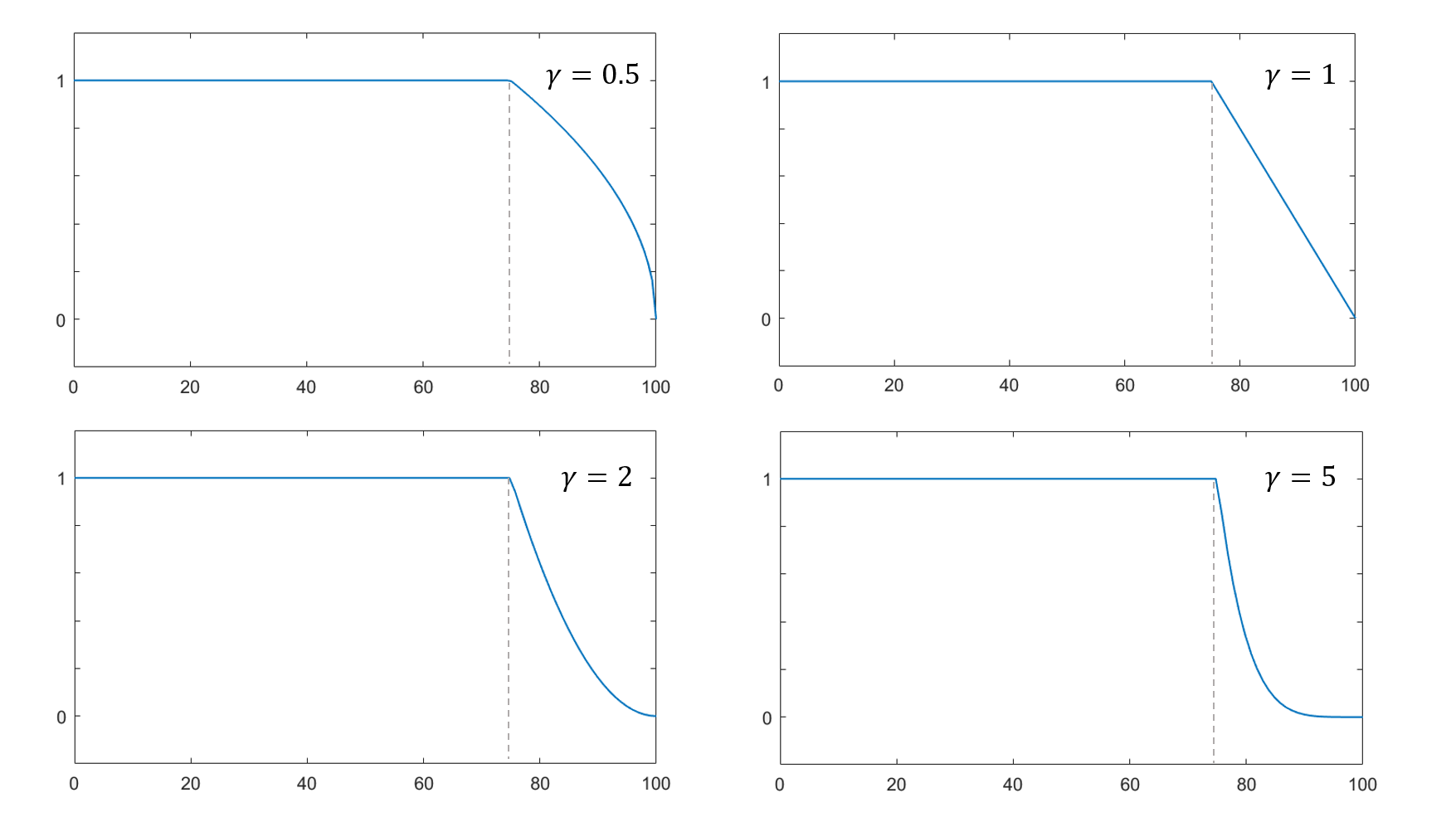}
    \caption{Plots of $f$ as a function of $\theta$ for a fixed $w > w_\text{vap}$ and varying values of $\gamma$.}
    \label{fig:f3}
\end{figure}

Thus, for mild temperatures $f(\theta,w) = 1$ and system \eqref{1} is retrieved. However, when $\theta$ reaches $\theta_{\text{vap}_1}$, $f$ decreases towards 0, damping the temperature growth. This process takes place until the point in which either the temperature has slowly increased up to $\theta_{\text{vap}_2}$ or almost all the water (except for a remaining $w_\text{vap}$) has vaporized. The exponent $\gamma$ only shapes the decrease in temperature growth, having a faster flattening for larger values of $\gamma$. Figure \ref{fig:f3} displays the plot of $f$ for fixed values of $\theta_{\text{vap}_1} = 75$ºC and $\theta_{\text{vap}_2} = 100$ºC and different exponents $\gamma$ for $w > \text{vap}$.

Note that this function $f$ only acts as a growth modulator for the energy arriving from the laser beam, accounted in the source $S$. When tissue water vaporization occurs, at 100ºC, the energy provided by the laser applicator does not contribute to a temperature increase (having $f(100$ºC$,w) = 0$ if $w > w_\text{vap}$, i.e. if the tissue water is still vaporizing). Instead, this energy is used for the phase change in the irradiated areas.

\subsection{Evolution of the tissue water concentration} \label{s1.3}

The remaining part of the model is devoted to the vaporization function $\Gamma$, modelling the loss of tissue water due to vaporization. Recall that the evolution equation for $w$ is of the form
\begin{equation} \label{5}
    \dfrac{\partial w}{\partial t}=-\Gamma \left( \theta \right) w,
\end{equation}
with $\Gamma(\theta)$ a suitable function. Integrating equation \eqref{5}, we obtain
$$
\displaystyle w(r,x,t) = w(r,x,0)\cdot \exp \left(-\int_0^t \Gamma(\theta(s)) ~ds\right).
$$
Therefore, for a fixed point $(r,t) \in \Omega_\text{cyl}$, the tissue water concentration suffers an exponential decay with time. Referring again to our previous argument, for mild temperatures the vaporization effects are almost negligible and almost no water is vaporized. As a result, $\Gamma$ should verify $\Gamma(\theta) \ll 1$ for low values of $\theta$ and then has to grow fast for $\theta$ between $\theta_{\text{vap}_1}$ and $\theta_{\text{vap}_2}$.

To obtain an explicit expression for $\Gamma(\theta)$ we make use of the experimental results presented in \cite{Yang-07}. Through measurements performed on \textit{ex-vivo} bovine liver, the authors fitted the following piecewise function $W(\theta)$ to represent the water content $W$ as a function of the temperature $\theta$.
\begin{equation} \label{6}
    W(\theta) = 778 \cdot \begin{cases}
        1 - \exp \left(\dfrac{\theta - 106}{3.42} \right) \quad & \text{if } \theta \leq 103\text{ºC}, \\\\
        0.03713 \theta^3 - 11.47 \theta^2 + 1182 \theta - 40582 \quad & \text{if } 103\text{ºC} < \theta \leq 104\text{ºC}, \\\\
        \exp \left( \dfrac{\theta - 80}{34.37} \right) \quad & \text{if } \theta > 104\text{ºC}.
    \end{cases}
\end{equation}
In our case, as we have tissue water concentration $w$ as a function of $r$, $x$ and $t$, through the chain rule, we obtain
$$
    \dfrac{d w}{d \theta} = \dfrac{\partial w}{\partial t} \cdot \dfrac{\partial t}{\partial \theta} = \dfrac{\partial w}{\partial t} \cdot \dfrac{1}{\partial \theta / \partial t}.
$$
We then substitute equations \eqref{3} on these partial derivatives. Moreover, as the experiments to obtain $W$ were performed \textit{ex-vivo}, we set the blood perfusion rate $\omega$ to zero and the source $S$ to a typical value $\Bar{S}$. We further assume that the diffusive term in the bio-heat equation for $\theta$ is small enough with respect to the other contributions, which will be justified on Section \ref{param} when realistic parameter values are considered. These simplifications yield
\begin{equation} \label{7}
    \dfrac{d w}{d \theta} = \dfrac{\partial w}{\partial t} \cdot \dfrac{1}{\partial \theta / \partial t} \approx \dfrac{-\Gamma \left( \theta \right) w,}{\Bar{S} \cdot f(\theta,w)}.
\end{equation}
Exploiting $W(\theta)$ given by \eqref{6}, we can solve for $\Gamma(\theta)$ in \eqref{7} and obtain
\begin{equation} \label{8}
    \Gamma(\theta) \approx - \dfrac{W'(\theta) \hspace{0.02 cm} \Bar{S} \hspace{0.02 cm} f(\theta, W(\theta))}{W(\theta)}.
\end{equation}
For the practical implementation of our model we use an exponential function which fits the approximation of $\Gamma(\theta)$ given by \eqref{8}.

\section{Parameters for porcine tissue} \label{param}

We adimensionalize system \eqref{3}. Denoting by $L$ and $R$ the length and radius of the considered tissue cylinder, we introduce the dimensionless temperature $\tilde{\theta}=\dfrac{\theta -\theta _{bl}}{\Delta \theta _{c}}$, with $\Delta \theta _{c}$ some maximal temperature difference with respect to blood, the dimensionless spatial coordinate $\tilde{x}=x/L$, $\tilde{r}=r/R$, the dimensionless time $\tilde{t}=t/t_{c}$, with $t_{c}$ a characteristic time, still to be selected, $\tilde{\Gamma}=\Gamma/\Gamma _{c}$ and $\tilde{S}=S/S _{c}$ with $\Gamma _{c}$ and $S_c$ a characteristic vaporization rate and incident irradiation intensity, still to be fixed as well.

With the aim of testing our model with currently existing experimental data regarding thermal ablation performed on porcine livers, we gather the parameter estimates for porcine liver tissue on Table \ref{t_estimates1}. 
\\
\begin{table}[h!]
\centering
\begin{tabular}{lll}
\toprule
\textbf{Parameter}                                   & \textbf{Value}                                         & \textbf{Source}   \\ \hline 
\rule{0pt}{3ex} $c_{ts}$, specific heat of the tissue                & 3640 J $\cdot$ kg$^{-1} \cdot$ K$^{-1}$                 & \cite{LITT}         \\
\rule{0pt}{3ex} $c_{bl}$, specific heat of blood                     & 3617 J $\cdot$ kg$^{-1} \cdot$ K $^{-1}$               & \cite{ITIS}   \\
\rule{0pt}{3ex} $\rho_{ts}$, mass density of the tissue             & 1137 kg $\cdot$ m$^{-3}$                               & \cite{LITT}        \\
\rule{0pt}{3ex} $\rho_{bl}$, mass density of blood                  & 1060 kg$\cdot$ m$^{-3}$                                & \cite{Fasano} \\
 \bottomrule
\end{tabular}
\caption{\label{t_estimates1} Parameter estimates for the densities and specific heats.}
\end{table}
\\
Next, we introduce a new dimensionless parameters $b$, along with two new time scales, $t_S$ and $t_D$, the characteristic diffusive time scale, as defined in Table \ref{t_adim}.
\begin{table}[h]
\centering
\begin{tabular}{ll}
\toprule
\textbf{Parameter} & \textbf{Expression}                                                                                                                            \\ \hline \rule{0pt}{4ex}  
$b$                & $\dfrac{\rho _{bl}c_{bl}}{\rho_{ts}c_{ts}}$                               \\  \rule{0pt}{5ex}  
$t_S$              & $\dfrac{\rho _{ts}c_{ts}\Delta \theta_{c}}{S_c }$      \\  \rule{0pt}{5ex}  
$t_D$              & $\dfrac{\rho_{ts} c_{ts} L^{2}}{k_{ts}}$     \\   
\bottomrule
\end{tabular}
\caption{\label{t_adim} New parameters for the model.}
\end{table}
\\
With these new parameters and variables, system \eqref{3} can by rewritten as \footnote{We have omitted the tildes to keep notation as light as possible.}
\begin{equation} \label{9}
\begin{cases}
    
\dfrac{\partial \theta }{\partial t}  =  \left( \dfrac{t_{c}}{t_{D}}
\right) \left[ \dfrac{\partial ^{2}\theta }{\partial x^{2}}+\dfrac{L^{2}}{
R^{2}}\dfrac{1}{r}\dfrac{\partial }{\partial r}\left( r\dfrac{\partial
\theta }{\partial r}\right) \right] - b\left( \dfrac{t_{c}}{1/\omega }\right)
\theta +\left( \dfrac{t_{c}}{t_{S}}\right) S \cdot f(\theta_{bl} + \Delta \theta_c \theta,w) ,  \\\\
\dfrac{\partial w}{\partial t}=-\left( \dfrac{t_{c}}{1/\Gamma _{c}}\right)
~\Gamma \left(\theta_{bl} + \Delta \theta_c \theta \right) w.
\end{cases}
\end{equation}
The system is considered for $(x,r)$ in the dimensionless version of $\Omega_{\text{cyl}}$ and for $t>0$. Notice that the Laplacian has been considered in cylindrical coordinates with no angular dependency due to the symmetry argument. 

System \eqref{9} highlights five distinct time scales: the initial characteristic time $t_c$, the characteristic diffusive time $t_{D}$ on the scale $L$, the perfusion time scale $1/\omega $, the evaporation time scale $1/\Gamma _{c}$ and the irradiation time scale $t_{S}$.

The remaining parameter values for the model not included on Table \ref{t_estimates1} can be found on Table \ref{t1}. The estimates of $k_{ts}$, $c_{ts}$ and $\rho_{ts}$ were obtained from \cite{LITT} and correspond to porcine liver tissue. Regarding $\omega$, the blood perfusion rate, its value varies considerably with the type of tissue and the species, \cite{Hemomath}. For pigs, in \cite{pigs} the mean blood perfusion rate in the liver is calculated for adult Göttingen minipigs and pre-pubertal Danish Landrace x Yorkshire pigs, obtaining 623/min and 950/min, respectively. As a preliminary approximation, we consider a mid value of $\omega = 750$ min$^{-1}$, though it would require a more concise estimate. For our domain, we take $L = 10$ cm and $R = 6$ cm. Lastly, for $\Delta \theta_c$, the maximal temperature difference, we take $80$ ºC, as our temperature scale will mainly measure between ambient temperature at the start of the experiment (around $20$ ºC) and vaporization temperatures reached at the end of the process (around $100$ ºC).

\begin{table}[h!]
\centering
\begin{tabular}{lll}
\toprule
\textbf{Parameter}                                   & \textbf{Value}                                         & \textbf{Source}   \\ \hline 
\rule{0pt}{3ex} $k_{ts}$, heat diffusivity coefficient in the tissue & 0.518 W $\cdot$ m$^{-1} \cdot$ K$^{-1}$                & \cite{LITT}     \\
\rule{0pt}{3ex} $\omega$, blood perfusion rate                        & 750 min$^{-1}$                          & \cite{pigs} \\
\rule{0pt}{3ex} $L$, cylinder length                       & $10 \cdot 10^{-2}$ m                          &  \\
\rule{0pt}{3ex} $R$, cylinder radius                       & $6 \cdot 10^{-2}$ m                      &  \\
\rule{0pt}{3ex} $\Delta \theta_c$,  maximal temperature difference                       & 80 ºK                      &  \\
\bottomrule
\end{tabular}
\caption{\label{t1} Remaining parameter estimates for the model.}
\end{table}

\subsection{Source function}

We are yet to set a source function $S$ modelling the laser irradiation. Though more accurate models can be selected, based on the \textit{radiative transfer equation}, (we refer the readers to \cite{9}), we select a simple Gaussian function $S$ based on an a quadratic exponential decay of the irradiation through the tissue.

Particularly, we take $S(x,r)$ as a product of two functions, $S_r(r)$ and $S_x(x)$, both representing said exponential decay in the $x$ and $r$ axis. To do so, we denote by $r_{\text{app}}$ the applicator radius and by $x_{\text{rad}_1}$ and $x_{\text{rad}_2}$ the two $x$ coordinates that bound the radiating part of the applicator. Thus, $\Gamma_{\text{rad}} = \{(x,r)\in \mathbb{R}^2: ~ x_{\text{rad}_1} < x < x_{\text{rad}_2}, ~ r = r_{\text{app}} \} \cup \{(x,r)\in \mathbb{R}^2: ~x =  x_{\text{rad}_2}, ~ 0 < r < r_{\text{app}}\}$. This way, we take
\begin{equation} \label{10}
 \displaystyle   S_r(r) = e^{-\beta (r-r_\text{app})^2},
\end{equation}
\begin{equation} \label{11}
    S_x(x) = \begin{cases}
        \displaystyle e^{-\beta (x-x_{\text{rad}_1})^2} \quad & \text{if } x \leq x_{\text{rad}_1}, \\\\
        1 \quad & \text{if } x_{\text{rad}_1} < x < x_{\text{rad}_2}, \\\\
        e^{-\beta (x-x_{\text{rad}_2})^2} \quad & \text{if } x_{\text{rad}_2} \leq x,
    \end{cases}
\end{equation}
where $\beta>0$ is a parameter representing the irradiation decay in the tissue. Lastly, we set
\begin{equation} \label{12}
    S(x,r) = S_0 \cdot S_x(x) \cdot S_r(r),
\end{equation}
with $S_0>0$ an irradiation intensity parameter that scales the function $S$, obtained through the applied laser power. Notice that when taking the dimensionless version of our model, we rescaled $S$ with $\Tilde{S} = S/S_c$ with $S_c$ a characteristic value of $S$. The most logical choice for $S_c$ is to simply take $S_c = S_0$ and then have $\Tilde{S}$ bounded between 0 and 1.

To obtain an initial estimate for $S_0$, given that for $S(x,r)$ we are considering a Gaussian beam, we have
$$S_0 = \frac{2P_0}{\pi a^2},$$
where $P_0$ is the laser power and $a$ is the beam waist. The typical laser power used for ablation procedures is between 20 and 40 W, while a preliminary estimate for the beam waist is $a = \mathcal{O}(1 \text{mm})$. This yields $S_0 = \mathcal{O}(10^7$ W $\cdot$ m$^{-2})$. The value of $\beta$ is later fitted to match the experimental data.

\subsection{Parameter fitting for the vaporization function}

\begin{figure}[htp]
    \centering
    \includegraphics[scale=0.6]{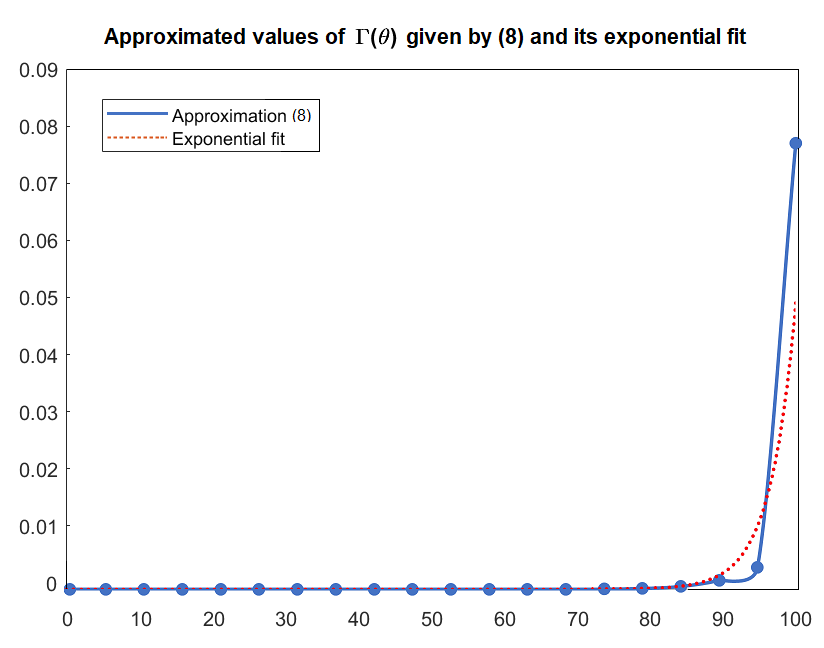}
    \caption{Approximation of $\Gamma(\theta)$ using \eqref{8} (in blue) and an exponential fit (dotted, in red).}
    \label{fig:f4}
\end{figure}

As a last step, we justify the assumptions made in \eqref{7} and fit a simpler function to our vaporization function $\Gamma$ approximated in \eqref{8}. The main feature was to approximate $S$ by a typical value $\Bar{S}$ and neglect the Laplacian, hence allowing us to obtain a function only dependant of $\theta$. The value of $\Bar{S}$ shall be taken just as $S_c$, as the irradiation intensity $S_0$. Regarding the Laplacian, we study the dimensionless version of the model, given by equations \eqref{9}. Recall that we had set $\omega = 0$ for \textit{ex-vivo} experiments. Hence, the coefficients of interest are $1/t_D$ and $1/t_S$, which can be estimated through the known parameters from Tables \ref{t_estimates1} and \ref{t1} and the previous result  $S_0 = \mathcal{O}(10^7$ W $\cdot$ m$^{-2})$. We find a significant enough difference of orders
$$\frac{1}{t_D} = \mathcal{O}(10^{-5} \text{ s}^{-1}), \quad \frac{1}{t_S} = \mathcal{O}(10^{-3} \text{ s}^{-1}),$$
which allows us to neglect the Laplacian in the approximation of $\Gamma$ given in \eqref{8}. Actually, to seek for a simpler alternative, we fit the approximated expression of $\Gamma(\theta)$ given by \eqref{8} with an exponential function. The plots of both functions are displayed together in Figure \ref{fig:f4}.
\\The obtained fit approximates expression \eqref{8} reasonably well, with a value of $R^2$ of $0.97$. For practical reasons we use the exponential fitting for the numerical simulations of the model, instead of the approximation \eqref{8}. As a final remark, for the value of $\Gamma_c$, though it is not entirely necessary, it can be taken as $10^{-2}$ to rescale the function.

\section{Numerical results} \label{s3}
Our study aims at comparing the model results with the experimental measurements from \cite{4}. As our main focus is on the vaporization effects, we compare our simulations with the case P34F47 in their study. For this case, a laser power of 34 W was employed during 1200 s, starting at room temperature of 20ºC. The experiments were performed on \textit{ex-vivo} porcine livers from previously slaughtered pigs, placing a temperature probe to measure the temperature evolution. In case P34F47, the probe position was $x =$ 23.8 mm, $r =$ 11.2 mm.

To compute the numerical solution to the dimensionless model \eqref{9}, a finite difference scheme was used, adapting the method to the geometry of $\Omega_{\text{cyl}}$ dividing it into two subrectangles. The dimensionless spatial step used was the same for $x$ and $r$, being $1/50$. The chosen time step was $4 \cdot 10^{-4}$ in turn. 

Given the \textit{ex-vivo} nature of the experiment, we set the blood perfusion rate $\omega = 0$. Regarding the remaining parameters, for the source function $\beta$ was fitted at $\beta = 40$ m$^{-2}$ and for the modulating function $f$, $\theta_{\text{vap}_1} = $ 75 ºC, $\theta_{\text{vap}_2} = $ 100 ºC, $w_\text{vap} = $5\% and $\gamma = 5$. The value of $S_0$ is approximated taking the slope of the initial measured temperature growth, making use of the smaller time scale of the diffusion. This yields $S_0 =$ 7.63 $\cdot$ 10$^{5}$. 

The initial conditions are taken spatially homogeneous, with initial temperature $\theta(x,r,0) = 20$ ºC and the tissue water concentration is set to $w(x,r,0) = W($20ºC$)\% =$ 77.8 \% in the entire domain. Notice that we are rescaling the tissue water concentration as a percentage.

\begin{figure}[htp]
    \centering
    \includegraphics[scale=0.5]{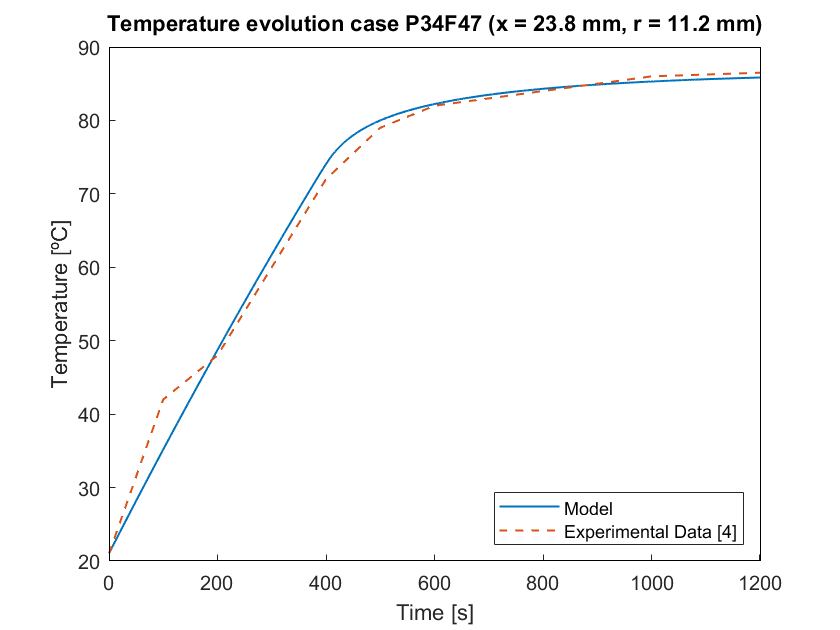}
    \caption{Numerical simulation of the temperature in case P34F47 plotted with the experimental measurements from \cite{4}.}
    \label{fig:f5}
\end{figure}

In Figure \ref{fig:f5} we plot in blue the obtained results for the temperature evolution at the probe position, together with the experimental measurements in the red dotted line.

\begin{figure}[htp]
    \centering
    \includegraphics[scale=0.465]{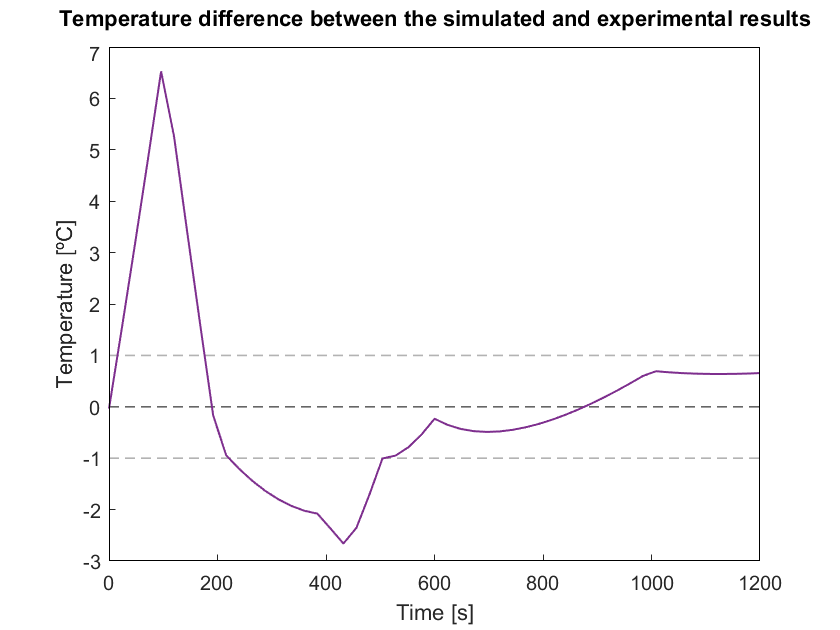}
    \caption{Obtained differences between the model simulations and the experimental results.}
    \label{fig:f6}
\end{figure}

As can be seen on the figure, the model seems to capture the temperature evolution well, especially the flattening due to the tissue water vaporization, occurring after approximately 400 seconds. After this change in the dynamics, the model fits the experimental results properly, with less than 1 ºC of difference. These differences are plotted in Figure \ref{fig:f6}. For the majority of the experiment, the difference between our predicted temperature and the measured data is of less than 1ºC, except on the early stages, before 200 seconds, where the measured temperature experiments a small bump, as well as a slight difference at the beginning of the vaporization effects.

\begin{figure}[htp]
    \centering
    \includegraphics[scale=0.7]{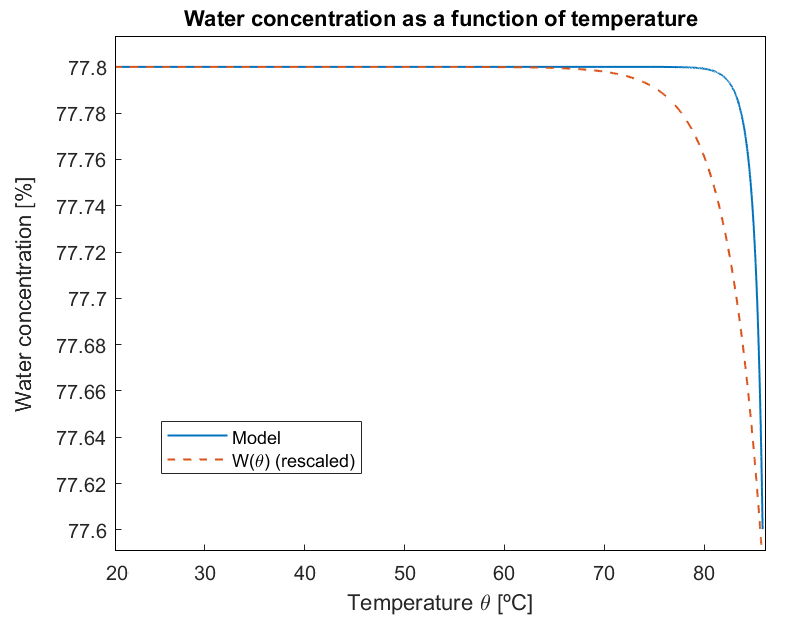}
    \caption{Numerical simulation of the temperature in case P34F47 plotted with the experimental measurements from \cite{4}.}
    \label{fig:f7}
\end{figure}

Lastly, we also plot the evolution of the tissue water concentration as a function of the temperature, to compare it with the graph of the experimentally fitted $W(\theta)$ from \eqref{6}. The results are depicted in Figure \ref{fig:f7}. The plots here reflect larger discrepancies, especially for temperatures above 70 ºC, where tissue water starts to vaporize in the model simulations. This however does not occur in the measurements before 80ºC. This can be a result of the approximation performed in \eqref{8}, particularly when approximating $S$ by $\Bar{S}$, to obtain a vaporization function that does not depend on the temperature. However, it could be reasonable to have a temperature-dependant vaporization rate, as certain areas are more exposed to the laser irradiance than others. For further improvements of the model, condensation effects can also be taken into account, following for instance the approach illustrated in \cite{LITT}. 

\section{Conclusion} \label{s4}
Throughout this work, a new model has been studied for describing laser-induced thermal ablation. The most important feature of the model is to consider a coupled system of differential equations describing not only the tissue temperature, but also its water concentration. In other models, like the one presented in \cite{LITT}, tissue water was taken as a known function, fitted from experimental data, whereas here it is taken as variable of the model. This allows for an alternative way of modelling the flatenning of the temperature growth as it approaches 100ºC, when the water starts to vaporize.

This way, we obtain system \eqref{3}, with the modulating function $f$ given in \eqref{4}. The loss of tissue water is represented by the vaporization function $\Gamma$, approximated in \eqref{8} using known experimental data.

The simulations are performed with realistic parameter values, collected on Tables  \ref{t_estimates1} and \ref{t1}. The obtained temperature is plotted (together with the measured data) in Figure \ref{fig:f5}. The fit is very good, as it can be seen in Figure \ref{fig:f6} where the differences are plotted. Lastly, the water concentration results as a function of the temperature were collected in Figure \ref{fig:f7}.

Overall, the model shows very good results when compared to the existing experimental data, except for minor discrepancies, mainly in the water content. More detailed terms can be incorporated as extensions of this model, but the proposed model seems to capture the physical reality reasonably well.

\end{document}